\documentclass[draft,english]{ourlematema}
\usepackage{amsmath}
\usepackage{verbatim}
\usepackage{mathtools}
\usepackage[makeroom]{cancel}
\usepackage{tikz}

\title{Points on rational normal curves and the ABCT variety}
\titlemark{Points on rational normal curves and the ABCT variety}

\titlenote{D.A. was supported by the DFG via the SFB-TRR 195 and the grant 530132094}

\MSC{05E05,14J81, 14N05, 14N15, 14N20}
\keywords{rational normal curve, configurations, Grassmannian, positive geometry, spinor helicity, scattering equations.}

\author{Daniele Agostini}
\address{%
Fachbereich Mathematik\\
Universit\"at T\"ubingen\\
\email{daniele.agostini@uni-tuebingen.de}
}
\author{Lakshmi Ramesh}
\address{%
Fakult\"at f\"ur Mathematik\\
Universit\"at Bielefeld\\
\email{lramesh@math.uni-bielefeld.de}
}
\author{Dawei Shen}
\address{%
Department of Mathematics\\
University of Michigan, Ann Arbor\\
\email{dwshen@umich.edu}
}
\date{2024/15/12}

\begin{document}

\maketitle

\begin{abstract}
The ABCT variety is defined as the closure of the image of $G(2,n)$ under the Veronese map. We realize the ABCT variety $V(3,n)$ as the determinantal variety of a vector bundle morphism. We use this to give a recursive formula for the fundamental class of $V(3,n)$. As an application, we show that  special Schubert coefficients of this class are given by Eulerian numbers, matching a formula by Cachazo-He-Yuan. On the way to this, we prove that the variety of configuration of points on a common divisor on a smooth variety is reduced and irreducible, generalizing a result of Caminata-Moon-Schaffler. 
\end{abstract}

\section{Introduction}

Let $k\leq n$ and let $d$ be such that $n\geq \binom{d+k-1}{d}$. The standard $d$-th Veronese embedding is the map $\mathbb{P}^{k-1}\to \mathbb{P}^{\binom{d+k-1}{d}-1}$ defined by $[x_1:x_2:\cdots:x_k]\mapsto [x_1^d:x_1^{d-1}x_2:\cdots :x_k^{d}]$, where the image point has coordinates running over all $d$ monomials of degree $d$ in $x_1,x_2,\cdots,x_k$. The $d$-th Veronese map on the Grassmannian $G(k,n)$ is the rational map $\theta_d\colon G(k,n) \dashrightarrow G\left( \binom{d+k-1}{d},n \right)$ given by
\begin{equation}\label{eq:veronesegrassmannianintro} 
\theta_{d}\colon 
\begin{pmatrix} 
x_{11} & x_{12} & \dots & x_{1n} \\
x_{21} & x_{22} & \dots & x_{2n} \\
\vdots & \vdots & \ddots & \vdots \\
x_{k1} & x_{k2} & \dots & x_{kn}
\end{pmatrix} \mapsto 
\begin{pmatrix} 
x_{11}^d & x_{12}^d & \dots & x_{1n}^d \\
x_{11}^{d-1} x_{21}  & x_{12}^{d-1}x_{22} & \dots & x_{1n}^{d-1}x_{2n} \\
\vdots & \vdots & \ddots & \vdots \\
x_{k1}^d & x_{k2}^{d} & \dots & x_{kn}^{d}
\end{pmatrix}.
\end{equation}
In other words, on each column of the matrix representing a point in $G(k,n)$, the map $\theta_d$ acts as the $d$-th standard Veronese embedding. Here, points in the two Grassmannians are given by the row spans of the corresponding matrices. Notice that this map is not defined precisely on the locus $Z_d(k,n)$ of matrices whose column vectors lie on a common hypersurface of degree $d$ in $\mathbb{C}^k$, see Remark \ref{rem:Zdnkcommonhyp}.

The \emph{ABCT variety} is the subvariety $V(d+1,n)\subseteq G(d+1,n)$ defined as the closure of the image of $\theta_{d}\colon G(2,n) \dashrightarrow G(d+1,n)$. 
\[ V(d+1,n) = \overline{\theta_{d}(G(2,n))} \]

This was introduced by Arkani-Hamed, Bourjaily, Cachazo and Trnka \cite{AHBCT11}, in relation to the spinor-helicity formalism in physics. See also \cite{lam2024moduli} and  \cite{EMPS24} for a recent discussion from a mathematical perspective.
From an algebro-geometric point of view, a general point in the variety $V(d+1,n)$ can be represented by matrices whose column vectors lie on a common rational normal curve in $\mathbb{P}^{d}$. Configuration of points on a rational normal curve are a classical topic, see for example \cite{W15}, but they have received renewed attention in recent years, see \cite{CGMS18},\cite{CGMS21},\cite{CMS23},\cite{CS21}, thanks also to connections with phylogenetics and data science. However, looking at these configurations inside the Grassmannian seems to be a novel perspective coming from physics. 

The purpose of this note is to investigate the geometry of the ABCT varieties, in particular we want to address some questions arising from Thomas Lam's lectures \cite{lam2024moduli}: what is the degree of $V(d+1,n)$ with respect to the Pl\"ucker embedding? More generally, what is the cohomology class of $V(d+1,n)$ in $G(d+1,n)$? What are its equations in Pl\"ucker coordinates? Is $V(d+1,n)$ a positive geometry?

Here we focus on the first significant case, that of the variety $V(3,n)$. We give an iterative formula to compute its cohomology class. Recall that the Schur polynomials $s_{\lambda_1,\lambda_2,\lambda_3}$ form a basis of the ring $\Lambda_3 =  \mathbb{Q}[x_1,x_2,x_3]^{\mathfrak{S}_3}$ of symmetric functions in three variables (see for example \cite[Corollary 1.1.4]{sturmfels2008algorithms}). There is a  surjective ring homomorphism $\Lambda_3 \to H^{*}(G(3,n)), f\mapsto [f]$ which sends the polynomial $s_{\lambda_1,\lambda_2,\lambda_3}$ to the Schubert class $[s_{\lambda_1,\lambda_2,\lambda_3}]$. We sometimes omit writing parts equal to $0$. For example, $s_{1,1}$ denotes $s_{1,1,0}$. Then the class of $V(3,n)$ is given as follows:

\begin{thm}\label{thm:classV3n}
Define recursively  the symmetric functions $f_m\in \Lambda_3$ as
\begin{align*}  
f_0 &=1, \quad f_1 = 4s_{1}, \quad f_2 = 11s_{2} + 6s_{1,1}, \\
f_m &= 2s_1\cdot f_{m-1} - (s_2+2s_{1,1})\cdot f_{m-2} + s_{2,1}\cdot f_{m-3} + 2^m \cdot f_{m}, \qquad \text{ for } m\geq 3.
\end{align*}
If $n\geq 5$,  the class of $V(3,n)$ in $G(3,n)$ is $[V(3,n)] = [f_{n-5}] \in H^{2n-10}(G(3,n))$.
\end{thm}
This recursion can be implemented easily on a computer algebra system, and it allows to compute explicit formulas in many cases: see Example \ref{ex:classesV3n}. Furthermore, from the expression of the class, it is straightforward to obtain the degree of $V(3,n)$ with respect to the Pl\"ucker embedding of $G(3,n)$: see Corollary \ref{cor:degreesV3n} and Example \ref{ex:degreeV3n}. We can also show in Corollary \ref{cor:euler} that the coefficient of $[s_{n-5}]$ in the expression of $[V(3,n)]$ is given  by the Eulerian number $A(n-3,1)$, matching a result of Cachazo, He and Yuan \cite{CHY13} (see also \cite[Theorem 6.4]{EMPS24} and  \cite[Section 4.4]{lam2024moduli}) in the spinor-helicity formalism. We elaborate on the connection to physics in Remark \ref{rem:physics} and Proposition \ref{prop:eulergeneral}.

The key to prove Theorem \ref{thm:classV3n} is to interpret $V(3,n)$ as a determinantal subvariety of $G(3,n)$. More precisely, we see that points in $V(3,n)$ can be represented by matrices whose columns lie on a degree two hypersurface in $\mathbb{C}^3$. This is precisely the locus $Z_2(3,n)$ where the map $\theta_2\colon G(3,n) \dashrightarrow G(6,n)$ of \eqref{eq:veronesegrassmannianintro} is not well-defined, which is the locus where the matrix defining $\theta_2$ has rank at most 5. In terms of the geometry of $G(3,n)$ this means that $Z_2(3,n)$ is the determinantal locus where a  morphism of vector bundles 
\[ \mathbb{C}^n \otimes \mathcal{O} \longrightarrow S^2\mathcal{U}^{\vee} \]	
has rank at most $5$. Here $\mathcal{U}$ is the universal subbundle on $G(3,n)$. Furthermore, we also show that this determinantal locus is also of expected codimension, so that its class can be computed via Porteous' formula, and we can then obtain Theorem \ref{thm:classV3n}. This gives the class of $Z_2(3,n)$ and to get the class of $V(3,n)$ we need to show that these coincide as schemes, meaning that $Z_2(3,n)$ is irreducible and reduced. In Theorem \ref{thm:Zdkn}, we prove this in general for any $Z_d(k,n) \subseteq G(k,n)$. The key is to use the main result of Caminata, Moon and Schaffler in \cite{CMS23}: they show that the variety of configurations of points in $\mathbb{P}^{k-1}$ lying in a common hypersurface of degree $d$ is irreducible and reduced. Here we need the following more general version of their statement:

\begin{thm}\label{thm:Xn}
Let $X$ be a smooth quasiprojective variety and let $V$ be a linear system on $X$ of dimension $\dim V = \ell$ such that the general divisor in $\mathbb{P}(V)$ is reduced and irreducible. For $n\geq 1$ define
\[ X^n(V) = \{ (p_1,\dots,p_n) \in X^n \,|\, p_1,\dots,p_n \text{ lie on a common divisor } D\in \mathbb{P}(V) \}. \]
If $n\geq \ell$, then $X^n(V)$ with its natural scheme structure is a reduced, Cohen-Macaulay and irreducible variety of expected codimension $n-\ell+1$ in $X^n$.
\end{thm}

In terms of commutative algebra, we get the following consequence:


\begin{cor}\label{cor:commalg}
Let $F_1,\dots,F_{\ell} \in \mathbb{C}[x_1,\dots,x_k]$ be linearly independent and assume that $\lambda_1 F_1 + \dots + \lambda_{\ell} F_{\ell}$
is irreducible for general $\lambda_1,\dots,\lambda_{\ell}\in \mathbb{C}$. If we choose $n \geq \ell$ different sets of variables $\mathbf{x}_j = (x_{1j},\dots,x_{kj})$ and we build the matrix
\[
\begin{pmatrix}
	F_1(\mathbf{x}_1) & F_1(\mathbf{x_2}) & \dots & F_1(\mathbf{x}_n) \\
	F_2(\mathbf{x}_1) & F_2(\mathbf{x_2}) & \dots & F_2(\mathbf{x}_n) \\
	\vdots & \vdots & \ddots & \vdots \\
	F_{\ell}(\mathbf{x}_1) & F_{\ell}(\mathbf{x_2}) & \dots & F_{\ell}(\mathbf{x}_n) \\
\end{pmatrix},
\]
then the ideal $I\subseteq \mathbb{C}[x_{ij}]$ of its maximal minors is prime of codimension $n-\ell+1$.
\end{cor}

The main result of \cite{CMS23} proves Theorem \ref{thm:Xn} for $X=\mathbb{P}^{k-1}$ and $V$ the complete linear system of forms of degree $d$. Equivalently, they prove Corollary \ref{cor:commalg} when $F_1,\dots,F_{\ell}$ are all the monomials of degree $d$ in $k$ variables: in this case, the matrix of Corollary \ref{cor:commalg} is precisely the matrix defining the map $\theta_d$ of \eqref{eq:veronesegrassmannianintro}, which is denoted by $\mathbf{M}_{r,d,n}$ in \cite{CMS23}. Our results, while apparently more general than those in \cite{CMS23}, can be proved in essentially the same way, but we felt that their statements might be worth writing down explicitly.

The variety $V(d+1,n)$ also invites interest as it is conjectured to be a positive geometry. We attempt to address Thomas Lam's conjecture in \cite{lam2024moduli} by considering the map $\theta_{d}: G(2,n) \dashrightarrow G(d+1,n)$, and understanding the face structure of the image that is inherited from $G(2,n)$. We consider the positroid stratification of $G(2,n)$ and analyse the subvarieties of $V(d+1,n)$ defined by the images of the positroid cells, which turn out to give us only a coarse stratification. We perform computations to show that dimensions can drop when we apply $\theta_{d}$ to positroid varieties in $ G(2,n)$, indicating that there might be boundaries that do not arise as the image of a positroid cell of $G(2,n)$.

The paper is structured as follows. In Section \ref{sec:configurations}, we prove Theorem \ref{thm:Xn} and Corollary \ref{cor:commalg}, which we then use in Section \ref{sec:vergrass} to prove in Theorem \ref{thm:Zdkn} that the determinantal loci $Z_d(k,n)$ are irreducible and reduced. In Section \ref{sec:V3n}, we prove Theorem \ref{thm:classV3n} and its corollaries. In Section \ref{sec:strat}, we conclude with discussions on whether $V(d+1,n)$ is a positive geometry.

\section{Configurations of points on a common divisor}\label{sec:configurations}

Let $X$ be a smooth quasiprojective complex variety of dimension $k$, with a line bundle $L$ and let $V\subseteq H^0(X,L)$ be a linear system with $\dim V = \ell$. If $\sigma\in V$ is a nonzero section, the vanishing locus $D=\{\sigma=0\}$ is a divisor on $X$. The points in $\mathbb{P}(V)$ can then be interpreted as divisors $D\subseteq X$. We will make the key assumption that \emph{a general divisor in $\mathbb{P}(V)$ is reduced and irreducible}. By Bertini's theorem, this is for example true if the base locus of $V$ has codimension at least $2$ and the image $\phi_V\colon X\dashrightarrow \mathbb{P}(V^{\vee})$ is not a curve. We consider the configuration of points lying on a common divisor in $\mathbb{P}(V)$.

\begin{dfn}[Configurations of points on a common divisor]\label{def:Xn}
	For $n\geq 1$ define
	\[ X^n(V) = \{ (p_1,\dots,p_n) \in X^n \,|\, p_1,\dots,p_n \text{ lie on a common divisor } D\in \mathbb{P}(V) \}. \]
\end{dfn} 

Since passing through a point is one linear condition on $\mathbb{P}(V)$,  we see that $X^n(V)=X^n$ whenever $n<\ell$. In the following, we will then consider the case where $n\geq \ell$. We consider the evaluation map, which is the following morphism of vector bundles
\begin{equation}\label{eq:evmap}  
\operatorname{ev}\colon V\otimes \mathcal{O}_{X^n} \longrightarrow \bigoplus_{i=1}^n \operatorname{pr}_i^*L, \qquad \sigma \mapsto (\sigma(p_1),\dots,\sigma(p_n)),  
\end{equation}
where $\operatorname{pr}_i\colon X^n \to X$ is the $i$-th projection.
We see that $(p_1,\dots,p_n)\in X^n(V)$ if and only if the evaluation map $\operatorname{ev}_{|(p_1,\dots,p_n)}$ has nonzero kernel, or, equivalently, if it does not have maximal rank. This means that $X^n(V)$ is the determinantal variety:
\begin{equation}\label{eq:Xndeglocus} 
X^n(V) = \left\{ (p_1,\dots,p_n) \in X^n \,|\, \operatorname{rk} \left(\operatorname{ev}_{|(p_1,\dots,p_n)}\right) < \ell  \right\}. 
\end{equation}
This endows $X^n(V)$ with a natural scheme structure: for $(p_1,\dots,p_n)\in X^n$ fix local coordinates $\mathbf{x}_j=(x_{1j},\dots,x_{kj})$ around each point $p_j$ (here we use that $X$ is smooth).  Fix also a basis $\sigma_1,\dots,\sigma_{\ell}$ of $V$. Around each of the $p_j$, the line bundle $L$ can be trivialized so that any section in $V$ can be written as an analytic function in the local coordinates: $\sigma_i = \sigma(x_{1j},\dots,x_{kj}) = \sigma(\textbf{x}_{j})$. Then the evaluation map in \eqref{eq:evmap} can be represented by a $n\times \ell$ matrix
\begin{equation}\label{eq:evmatrix}
\operatorname{ev} = 
\begin{pmatrix} 
\sigma_1(\mathbf{x}_1) & \sigma_2(\mathbf{x}_1) & \dots & \sigma_\ell(\mathbf{x}_1) \\
\sigma_1(\mathbf{x}_2) & \sigma_2(\mathbf{x}_2) & \dots & \sigma_\ell(\mathbf{x}_2) \\
\vdots & \vdots & \ddots & \vdots \\
\sigma_1(\mathbf{x}_n) & \sigma_2(\mathbf{x}_n) & \dots & \sigma_\ell(\mathbf{x}_n)
\end{pmatrix},
\end{equation} 
so that $X^n(V)$ is locally defined by the vanishing of maximal $\ell\times \ell$ minors of this matrix. We will consider $X^n(V)$ with this scheme structure. We start by showing that $X^n(V)$ is irreducible of the expected codimension. To do so, it is natural to consider the incidence correspondence
\begin{align} 
\Sigma^n(V) &= \{ (D,(p_1,\dots,p_n)) \in \mathbb{P}(V) \times X^n \,|\, p_j\in  D \quad \text{ for } j=1,\dots,n \}. 
\end{align}
If we fix a basis $\sigma_1,\dots,\sigma_{\ell}$ of $V$, then each element in $\mathbb{P}(V)$ is a linear combination $\sum_{h=1}^{\ell} a_h \sigma_h$ up to a scalar, and $\Sigma^n(V)$ is defined on $\mathbb{P}(V)$ by the conditions $\left( \sum_{h=1}^{\ell} a_h \sigma_h \right)(p_j)=0$ for $j=1,\dots,n$. Each of these conditions defines a divisor of class $\operatorname{pr}_{\mathbb{P}(V)}^*\mathcal{O}_{\mathbb{P}(V)}(1) \otimes \operatorname{pr}_j^*L$ on $\mathbb{P}(V)\times X^n$. Hence $\Sigma^n(V)$ is an intersection of $n$ divisors and acquires a natural scheme structure in this way. More concretely, if $\mathbf{x}_j$ are local coordinates around the points $p_j$ as before, we can write the conditions defining $\Sigma^n(V)$ as
\begin{equation}\label{eq:evmatrixzero}
\begin{pmatrix} 
\sigma_1(\mathbf{x}_1) & \sigma_2(\mathbf{x}_1) & \dots & \sigma_\ell(\mathbf{x}_1) \\
\sigma_1(\mathbf{x}_2) & \sigma_2(\mathbf{x}_2) & \dots & \sigma_\ell(\mathbf{x}_2) \\
\vdots & \vdots & \ddots & \vdots \\
\sigma_1(\mathbf{x}_n) & \sigma_2(\mathbf{x}_n) & \dots & \sigma_\ell(\mathbf{x}_n)
\end{pmatrix}
\begin{pmatrix}
a_{1}\\ a_2 \\ \vdots\\ a_{\ell}
\end{pmatrix}=0.
\end{equation}

\begin{prop}\label{prop:Sigman}
	Let $X$ be a smooth quasiprojective variety and let $V$ be a linear system on $X$ of dimension $\dim V = \ell$ such that the general divisor in $\mathbb{P}(V)$ is reduced and irreducible. Assume that $n\geq \ell$ and consider the incidence correspondence $\Sigma^n(V)$ with the scheme structure mentioned above. 
	\begin{itemize}
		\item[(i)] $\Sigma^n(V)$ is a reduced and irreducible complete intersection of expected codimension $n$.
		\item[(ii)] A point $(D,(p_1,\dots,p_n))$ is smooth on $\Sigma^n(V)$ if and only if the $p_i$ contained in the singular locus $D^{\operatorname{sing}}$
are pairwise distinct points, imposing independent conditions on the linear system $V$. 	
\end{itemize} 
\end{prop}
\begin{proof}
(i). For the irreducibility and the codimension, we follow a different approach than \cite[Lemma 3.2]{CMS23}. Since $\Sigma^n(V)$ is the intersection of $n$ divisors, each irreducible component has codimension at most $n$, meaning dimension at least $nk-n+\ell-1$. Consider the projection $\operatorname{pr}_{\mathbb{P}(V)}\colon \Sigma^n(V) \to \mathbb{P}(V)$: the fiber over a divisor $D\in \mathbb{P}(V)$ is given by $D^n$, hence it has dimension $nk-n$. Thus, every irreducible component of $\Sigma^n(V)$ has dimension at most $nk-n+\ell-1$. This proves that every component has dimension exactly $nk-n+\ell-1$. Furthermore, thanks to our assumptions on $V$, there is a nonempty open subset $U\subseteq \mathbb{P}(V)$ such that every $D\in U$ is  irreducible. The fibers of $\operatorname{pr}_{\mathbb{P}(V)}$ over $U$ are then all irreducible of the same dimension. Hence, the Zariski closure of $\operatorname{pr}_{\mathbb{P}(V)}^{-1}(U)$ in $\Sigma^n(V)$ is irreducible of dimension $nk-n+\ell-1$ and it must be a component. If there is another component $\Sigma'$, then the closure of $\operatorname{pr}_{\mathbb{P}(V)}(\Sigma')$ must be a proper subset of $\mathbb{P}(V)$, but then $\Sigma'$ would have dimension smaller than $nk-n+\ell-1$, a contradiction. This shows that $\Sigma^n(V)$ is irreducible and of expected codimension $n$, and since it is an intersection of $n$ divisors, it must be a complete intersection. The rest of the proof proceeds as in \cite[Lemma 3.3]{CMS23}: since it is a complete intersection,  it is Cohen-Macaulay, hence has no embedded components. Thus, to prove that it is reduced, it is enough to show that there is one smooth point. This follows from point (ii): it is enough to take one general divisor $D\in \mathbb{P}(V)$ and smooth points $p_1,\dots,p_n\in D$. Such points exist since $D$ is generically smooth.

(ii) This follows by the Jacobian criterion, using the equations in \eqref{eq:evmatrixzero}. The Jacobian matrix can be computed as in  \cite[Proof of Lemma 3.3]{CMS23}, and to conclude one can reason as in   \cite[Proof of Theorem 4.6]{CMS23}. 
\end{proof}
Now we can give the proofs of Theorem \ref{thm:Xn} and Corollary \ref{cor:commalg}.
\begin{proof}[Proof of Theorem \ref{thm:Xn}]
	The proof is as in \cite{CMS23}: the variety $X^n(V)$ is the image of the projection $\operatorname{pr}_{X^n}\colon \Sigma^n(V) \to X^n$, and then Proposition \ref{prop:Sigman} shows that $X^n(V)$ is irreducible and of dimension at most $nk-n+\ell-1$, hence of codimension at least $n-\ell+1$. However, since $X^n(V)$ is a degeneracy locus of a morphism from a bundle of rank $\ell$ to a bundle of rank $n$, it must have codimension at most $n-\ell+1$. This proves that it has exactly this codimension. In particular, it is Cohen-Macaulay. At this point, one can prove as in \cite[Lemma 3.4]{CMS23} that the map $\operatorname{pr}_{X^n}\colon \Sigma^n(V) \to X^n(V)$ is an isomorphism over the open set in $X^n(V)$ of the configurations $(p_1,\dots,p_n)$ contained in exactly one divisor in $\mathbb{P}(V)$. To adapt the proof of \cite{CMS23} to our setting,  it is enough to replace their matrix $\mathbf{M}_{r,d,n}$ with the matrix of \eqref{eq:evmatrix}, and their monomials of degree $d$ with the basis $\sigma_1,\dots,\sigma_{\ell}$ of $V$. In particular, $X^n(V)$ is generically reduced, and since it is Cohen-Macaulay, it is everywhere reduced. 
\end{proof}

\begin{proof}[Proof of Corollary \ref{cor:commalg}]
Consider the variety $X=\mathbb{C}^k$ and the linear system $V=\langle F_1,\dots,F_{\ell}\rangle$. This linear system fulfills the hypotheses of Theorem \ref{thm:Xn}, so that $X^n(V)$ is reduced and irreducible, meaning that the ideal defining it as a scheme is prime. This ideal is precisely the ideal $I$ of maximal minors of the matrix appearing in our statement.
\end{proof}

\section{Veronese maps on the Grassmannian}\label{sec:vergrass}
Let $k\leq n$ be two positive integers. The Grassmannian $G(k,n)$ parametrizes $k$-dimensional subspaces of a 
complex vector space of dimension $n$. We recall some basic facts about this space. First, we will represent a point in $G(k,n)$ via a complex $k\times n$ matrix of rank $k$
\begin{equation}\label{eq:matrixrepresentation}
M = 
\begin{pmatrix} 
x_{11} & x_{12} & \dots & x_{1n} \\
x_{21} & x_{22} & \dots & x_{2n} \\
\vdots & \vdots & \ddots & \vdots \\
x_{k1} & x_{k2} & \dots & x_{kn}
\end{pmatrix}
\end{equation}
so that $W = \operatorname{rowspan}(M)$ is a $k$-dimensional subspace of $V=\mathbb{C}^{1\times n}$. Two such matrices have the same row span if and only they differ by a left multiplication by an invertible $k\times k$ matrix. This realizes the Grassmannian as a quotient $G(k,n) = \operatorname{GL}_k\setminus \operatorname{Mat}_{k\times n}^{\circ}$, where $\operatorname{Mat}^{\circ}_{k\times n}$ is the set of $k\times n$ matrices of maximal rank. Sometimes we will also use the notation $[M]\in G(k,n)$ for a subspace represented by a matrix.

The Grassmannian comes equipped with a natural tautological bundle $\mathcal{U}$ of rank $k$ whose fiber at a point $W\in G(k,n)$ is $\mathcal{U}_{|W} = W$, the vector space itself. Correspondingly, there is also a rank $n-k$ quotient bundle $\mathcal{Q}$ with fiber $\mathcal{Q}_{|W} \cong V/W$. These fit into the Euler exact  sequence
\begin{equation}\label{eq:undualeulersequence}
 0 \longrightarrow \mathcal{U} \longrightarrow V\otimes \mathcal{O} \longrightarrow \mathcal{Q} \longrightarrow 0
\end{equation}
which globalizes the exact sequence of vector spaces $0\to W \to V \to V/W \to 0$. 

\begin{rem}\label{rem:matrixtheta}
The Euler sequence \eqref{eq:undualeulersequence} allows to interpret the matrix $M$ of Equation \eqref{eq:matrixrepresentation} as an object on $G(k,n)$. Dualizing the Euler sequence, we obtain another exact sequence
\begin{equation}\label{eq:eulersequence}
    0 \longrightarrow \mathcal{Q}^{\vee} \longrightarrow V^{\vee}\otimes \mathcal{O} \overset{\cdot M}{\longrightarrow} \mathcal{U}^{\vee} \longrightarrow 0
\end{equation}
The last map takes a globally defined functional $\phi\in V^{\vee}$ and, over each $W\in  G(k,n)$ restricts it to a subspace $W = \mathcal{U}_{|W}$. This map is represented by the matrix $M$ of Equation \eqref{eq:matrixrepresentation}: indeed, if $W = \operatorname{rowspan}(M)$, a basis of $\mathcal{U}_{|W} = W$ is given by the rows $M_1,\dots,M_k$ and there is a corresponding dual basis $M_1^{\vee},\dots,M_k^{\vee}$ of $\mathcal{U}^{\vee}_{|W} = W^{\vee}$.  A basis of $V^{\vee} = \mathbb{C}^n$ is given by the usual canonical column vectors $e_1,\dots,e_n$ and the coordinates of the restriction ${e_j}_{|W}$ with respect to the basis $M_1^{\vee},\dots,M_k^{\vee}$ are given precisely by the $j$-th column of $M$. 
\end{rem}
Fix now a positive integer $d$ and consider the map $\theta_{d}^{k,n}$ (we will often use only $\theta_d$ for simplicity ) defined on the $k\times n$ matrices of Equation \eqref{eq:matrixrepresentation}, as follows:
\begin{equation}\label{eq:veronesegrassmannian} 
\theta_{d}^{k,n}\colon M \mapsto \theta_d^{k,n}(M) = 
\begin{pmatrix} 
x_{11}^d & x_{12}^d & \dots & x_{1n}^d \\
x_{11}^{d-1} x_{21}  & x_{12}^{d-1}x_{22} & \dots & x_{1n}^{d-1}x_{2n} \\
\vdots & \vdots & \ddots & \vdots \\
x_{k1}^d & x_{k2}^{d} & \dots & x_{kn}^{d}
\end{pmatrix}.
\end{equation}
In other words, on each column of the matrix, the map $\theta_d$ acts as the $d$-th standard Veronese embedding.  Furthermore, if $A\in \operatorname{GL}_k$, then $\theta_d(A\cdot M) = S^dA \cdot \theta_d(M)$, where $S^dA \in \operatorname{GL}_{\binom{d+k-1}{k-1}}$ is a symmetric power of $A$. To see this, consider the usual Veronese embedding as $v\colon \mathbb{C}^k \to S^d\mathbb{C}^k, v\mapsto v^d$: then $S^dA$ is the matrix form of the of the usual action of $A\in \operatorname{GL}_k$ on $S^d\mathbb{C}^k$. One can then define the Veronese maps seen in the introduction.

\begin{dfn}[Veronese maps on the Grassmannian]\label{def:veronesegrassmannian}
	Assume that $\binom{d+k-1}{d}\leq n$. The $d$-th Veronese map on the Grassmannian $G(k,n)$ is the rational map
\[ \theta_d^{k,n} \colon G(k,n) \dashrightarrow G\left( \binom{d+k-1}{d},n \right), \qquad [M] \mapsto [\theta_d^{k,n}(M)] \]
\end{dfn}

\begin{rem}\label{rem:matrixtthetaMbundles}
	Following the same reasoning of Remark \ref{rem:matrixtheta}, we can also interpret the matrix $\theta_d(M)$ as a morphism of vector bundles:
	\begin{equation}\label{eq:morphismbundles} 
 V^{\vee}\otimes \mathcal{O} \xrightarrow{\cdot \theta_d(M)} S^d\mathcal{U}^{\vee} 
 \end{equation}
    The morphism is defined on the canonical basis vectors $e_1,\dots,e_n\in \mathbb{C}^n = V$ as follows: the restriction of $e_j$ to $W\in \operatorname{G}(k,n)$ defines a functional in $W^{\vee}$. Then the restriction of $e_j^d\in S^dV$ to $W$ defines a functional in $S^dW^{\vee} = S^d\mathcal{U}_{|W}$. To see that this morphism is represented by the matrix $\theta_d(M)$ assume that $W\in \operatorname{G}(k,n)$ is represented as $W = \operatorname{rowspan}(M)$, then a basis of $S^dW^{\vee}$ is given by the degree $d$ monomials $(M_1^{\vee})^d,(M_1^{\vee})^{d-1}(M_2),\dots,(M_k^{\vee})^d$, where the $M_i$ are the rows of $M$ as in Remark \ref{rem:matrixtheta}. The coordinates of $(e_j)^d$ with respect to the above basis of $S^dW^{\vee}$ are given precisely by the $j$-th column of $\theta_d(M)$.
\end{rem}

The Veronese map $\theta_d$ fails to be well-defined exactly on the locus where the matrix $\theta_d(M)$ does not have maximal rank. This is by definition also the degeneracy locus of the morphism \eqref{eq:morphismbundles}. We denote it as
\begin{equation}
Z_d(k,n) = \left\{ [M] \in G(k,n) \,|\, \operatorname{rk}\theta_{d}(M) < \binom{d+k-1}{d} \right\}.
\end{equation}
This has a scheme structure given by the vanishing of the maximal minors of the matrix $\theta_d(M)$ of \eqref{eq:veronesegrassmannian}.

\begin{rem}\label{rem:Zdnkcommonhyp}
On $\mathbb{C}^k$, consider the linear system $V_d = \mathbb{C}[x_1,\dots,x_k]_d$ of all forms of degree $d$, and take the configuration space $(\mathbb{C}^k)^n(V_d)$ as in Definition \ref{def:Xn}. Then $Z_d(k,n)$ is, as a set, the quotient of $(\mathbb{C}^k)^n(V_d) \cap \operatorname{Mat}^{\circ}(k,n)$ by $\operatorname{GL}_k$: this is because the maximal minors of $\theta_d(M)$ are the same as the maximal minors of $\theta_d(M)^t$. In other words, a point in $Z_d(k,n)$ is represented by a matrix $M$ whose columns lie on a common degree $d$ hypersurface.
\end{rem}

\begin{thm}\label{thm:Zdkn}
Assume  $d\geq 2$ and $n\geq \binom{d+k-1}{d}$. Then $Z_d(k,n)$ with its natural scheme structure is reduced, irreducible and Cohen-Macaulay of expected codimension $n-\binom{d+k-1}{d}+1$ in $G(k,n)$.
\end{thm}
\begin{proof}
We can restrict ourselves to prove the statement for the intersections of $Z_d(k,n)$ with the standard open charts of $G(k,n)$. By symmetry it is enough to consider the open chart  $\mathbb{C}^{k(n-k)} \subseteq G(k,n)$ where the last $k$ columns of $M$ in \eqref{eq:matrixrepresentation} give the identity matrix. On this chart, the matrix $\theta_d(M)$ has the form:
\begin{equation*}
    \theta_d(M) = \begin{pmatrix} 
     x^d_{11} & x^d_{12} & \dots & x^d_{1,n-k} & 1 & \dots & 0 \\ 
     x^{d-1}_{11}x_{21} & x^{d-1}_{12}x_{22} & \dots & x^{d-1}_{1,n-k}x_{2,n-k} & 0 & \dots & 0 \\
     \vdots & \vdots & \ddots & \vdots & \vdots & \ddots & \vdots \\
    x^d_{k1} & x^d_{k2} & \dots & x^d_{k,n-k} & 0 & \dots & 1 \\ 
    \end{pmatrix}
\end{equation*}
To compute the ideal generated by the $k\times k$ minors, we can subtract the last $k$ columns from the previous one, so that all entries of the form $x_{ij}^d$ are zero. Then, the ideal generated by the maximal minors of the matrix $\theta_d(M)$ coincides with the ideal of the maximal minors of the matrix
\begin{equation}\label{eq:thetadprime}
\theta_d'(M) = \begin{pmatrix} x_{11}^{d-1}x_{21} & x_{12}^{d-1}x_{22} & \dots & x_{1,n-k}^{d-1}x_{2,n-k}\\
\vdots & \vdots & \ddots & \vdots \\
x_{k-1,1}x_{k1}^{d-1} & x_{k-1,2}x_{k2}^{d-1} & \dots & x_{k-1,n-k}x_{k,n-k}^{d-1}
\end{pmatrix}.
\end{equation} This is the matrix $\theta_d(M)$ but without the monomials of the form $x_{ij}^d$. Of course, the ideal of maximal minors of $\theta'_d(M)$ coincides with the ideal of maximal minors of $\theta'_d(M)^t$. We can phrase this in terms of the configuration spaces of Definition \ref{def:Xn} as follows: on $\mathbb{C}^{k}$ consider the linear system given by $V = \langle x_{1}^{d-1}x_{2},\dots, x_{k-1}x_k^{d-1} \rangle$. This consists of all degree $d$ forms passing through the canonical basis vectors in $\mathbb{C}^k$. The previous discussion shows that $Z_d(k,n) \cap \mathbb{C}^{k(n-k)} = (\mathbb{C}^{k})^{n-k}(V)$.
Our assumptions guarantee that $n-k\geq \dim V$, furthermore, it is straightforward to show that a general hypersurface in $\mathbb{P}(V)$ is smooth and irreducible. Then the result follows from Theorem \ref{thm:Xn}.
\end{proof}

\begin{rem}
If $n<\binom{d+k-1}{k-1}$, then $Z_d(k,n) = G(k,n)$. If instead $d=1$ then it follows from the definition of the Grassmannian that $Z_1(k,n)=\emptyset$.
\end{rem}

\begin{rem}
Since the degeneracy locus $Z_d(k,n)$ has the expected codimension, we can use Porteous' formula to compute its class $[Z_d(k,n)]$ in the Chow ring of $G(k,n)$. In the next section, we are going to apply this method to compute the class of the variety $V(3,n) \subseteq G(3,n)$. 
\end{rem}

We conclude this section with a consequence  in commutative algebra:

\begin{cor} Assume $d\geq 2$ and $n\geq \binom{d+k-1}{d}$. Then the ideal generated by the maximal minors of the matrix $\theta_d'(M)$ of Equation \eqref{eq:thetadprime} is prime.
\end{cor}
\begin{proof}
By the proof of Theorem \ref{thm:Zdkn}, we see that the maximal minors of $\theta_d'(M)$ define the ideal of the scheme $Z_d(k,n)\cap \mathbb{C}^{k(n-k)}$. Since this scheme is reduced and irreducible, its ideal must be prime. 
\end{proof}
We also give a purely algebraic proof that the ideal is radical for a case particularly interesting for us:  $k=3$ and $d=2$. We are considering the ideal of maximal minors of the matrix
\begin{equation}\label{matrixtheta'}
\theta_d'(M) = \begin{pmatrix}
x_{11}x_{21} & x_{12}x_{22} & \dots & x_{1,n-3}x_{2,n-3} \\
x_{11}x_{31} & x_{12}x_{32} & \dots & x_{1,n-3}x_{3,n-3} \\
x_{21}x_{31} & x_{22}x_{32} & \dots & x_{2,n-3}x_{3,n-3} 
\end{pmatrix}. 
\end{equation}

\begin{prop}
Let $I \subseteq \mathbb{C}[x_{ij}|i\in [3],j\in [n]]$ be the ideal of maximal minors of the matrix $\theta'_d(M)$ in (\text{\ref{matrixtheta'}}). Consider a lexicographic term order with
\[ x_{11}>x_{12}>\dots>x_{1,n-3}>x_{21}>x_{22}>\dots>x_{3,n-3} \]
Then the maximal minors of $\theta_d'(M)$ form a Gr\"obner basis for $I$ with respect to this term order. In particular, $I$ is radical.
\end{prop}
\begin{proof}
Let $f_{ijk}$ and $f_{rst}$ be two minors of the above matrix corresponding to the columns $1\leq i<j<k\leq n-3$ and $1\leq r<s<t\leq n-3$. We want to use Buchsberger's criterion and prove that the $S$-polynomial of these two minors reduces to zero under division by the other minors. If this holds in the ring $\mathbb{C}[x_{ab}|a\in[3],b\in \{i,j,k,r,s,t\}]$, then it holds in the larger ring as well, since the above term order restricts in a compatible way. Thus, it is enough to check the statement when the matrix $\theta'_d(M)$ has $6$ columns, and this can be checked quickly by OSCAR \cite{OSCAR} or Macaulay2 \cite{M2}.
Since each maximal minor is a sum of squarefree terms, it follows that the initial ideal of $I$ is radical, so $I$ is radical as well.
\end{proof}

\section{The variety $V(3,n)$}\label{sec:V3n}

Let $n\geq d+1$. The subvariety $V(d+1,n)\subseteq G(d+1,n)$ is the Zariski closure of the image of the Veronese map
\[ \theta_{d}^{2,n}\colon G(2,n) \dashrightarrow G(d+1,n), \qquad \begin{pmatrix} x_1 & \dots & x_n \\ y_1 & \dots & y_n \end{pmatrix} \mapsto \begin{pmatrix} x_1^{d} & \dots & x_n^{d} \\ x_1^{d-1}y_1 & \dots & x_n^{d-1}y_n \\\vdots & & \vdots \\ y_1^{d} & \dots & y_n^{d} \end{pmatrix} \]
as in (\ref{eq:veronesegrassmannian}). By construction, a general point in $V(d+1,n)$ is represented by a matrix $[M]\in G(d+1,n)$ whose column vectors lie on a common rational normal curve in $\mathbb{P}^{d}$. We will always take $V(d+1,n)$ with its reduced scheme structure. We focus on the case $V(3,n)$. For $3\leq n\leq 5$, we have $V(3,n)=G(3,n)$. A simple but important observation is:

\begin{prop}\label{prop:V=degeneracy}
Let $n\geq 3$. The varieties $V(3,n)$ and $Z_2(3,n)$ coincide in $G(3,n)$.
\end{prop}
\begin{proof}
Recall from \cite[Definition 4.8]{lam2024moduli} that $V(3,n)$ is irreducible of dimension $2n-4$. In particular, it is the whole of $G(3,n)$ if $n<6$, and the same is true of $Z_2(3,n)$. If instead $n\geq 6$, then we know from Theorem \ref{thm:Zdkn}, that $Z_{2}(3,n)$ is reduced and irreducible of dimension $2n-4$. Hence it is enough to prove that $V(3,n)\subseteq Z_2(3,n)$. To show this, observe that a general point in  $V(3,n)$ is represented by a matrix $[M] \in G(3,n)$ whose projectivized column vectors lie on a common rational normal curve in $\mathbb{P}^2$, see Remark \ref{rem:Zdnkcommonhyp}. Since such a curve is a plane conic, we see that $[M]\in Z_2(3,n)$.
\end{proof}

It is then straightforward to compute defining equations for $V(3,n)$:

\begin{cor}
    If $n\geq5$. Then $V(3,n)$ is cut out as a scheme by the following quartic equations in the Pl\"ucker coordinates on $G(3,n)$:
    \[p_{i_1,i_2,i_3}p_{i_1,i_5,i_6}p_{i_2,i_4,i_6}p_{i_3,i_4,i_5}-p_{i_2,i_3,i_4}p_{i_1,i_2,i_6}p_{i_1,i_3,i_5}p_{i_4,i_5,i_6} = 0,\] where $i_1<i_2<i_3<i_4<i_5<i_6$ vary across all $6$-element subsets of $\{1,\dots,n\}$.
\end{cor}
\begin{proof}
Since $V(3,n)=Z_2(3,n)$, it is cut out by the $6\times 6$ minors of the matrix defining $\theta^{3,n}_2$. One can check, for example with a computer algebra system, that these minors can be expressed in Pl\"ucker coordinates by the quartics above (see also \cite[Remark 3.3]{CGMS18}).
\end{proof}

We now proceed to compute the class of $V(3,n)$ and prove Theorem \ref{thm:classV3n}. Recall (see  \cite[Section 9.4]{fulton1997young}) that the cohomology ring $H^{*}(G(3,n))$ can be seen as a quotient of the ring of symmetric functions $\Lambda_3 = \mathbb{Q}[x_1,x_2,x_3]^{\mathfrak{S}_3}$ via a map $\Lambda_3 \to H^{*}(G(3,n)), f\mapsto [f]$ that sends the Schur polynomial $s_{\lambda}$ of a partition $\lambda=(\lambda_1,\lambda_2, \lambda_3)$ to the class $[s_{\lambda}]$ of any Schubert variety associated to $\lambda$. The $[s_{\lambda}]$ are classes of codimension $|\lambda| = \sum_i \lambda_i$ and $[s_{\lambda}] \ne 0$ if and only if $\lambda$ fits into a $3\times (n-3)$ box, meaning that  $n-3 \geq \lambda_1 \geq \lambda_2 \geq \lambda_3 \geq 0$.  

We start with a lemma on symmetric functions: for a partition $\lambda$, denote by $e_{\lambda},h_{\lambda}$ the elementary and  complete symmetric functions in three variables associated to  $\lambda$. Define also symmetric functions $f_d(x_1,x_2,x_3)$ via the generating series
\begin{equation}\label{eq:genseries} 
F(x_1,x_2,x_3,t)= \sum_{d=0}^{\infty} f_d(x_1,x_2,x_3) \cdot t^d = \prod_{i=1}^3 \frac{1}{1-2x_i t} \cdot \prod_{1\leq i < j \leq 3} \frac{1}{1-(x_i+x_j)t}
\end{equation}

\begin{lemma}\label{lem:symmfunctions}
With the previous notation, it holds that
\begin{align*}  
f_0 &= 1, \quad f_1 = 4s_{1}, \quad f_2 = 11s_{2} + 6s_{1,1} \\
f_m &= 2s_1\cdot f_{m-1} - (s_2+2s_{1,1})\cdot f_{m-2} + s_{2,1}\cdot f_{m-3} + 2^m \cdot s_{m}, \qquad \text{ for } m\geq 3.
\end{align*}
\end{lemma}

\begin{proof}
Recall that $\prod_{i=1}^3(1-x_i t)^{-1} = \sum_{m=0}^{\infty} h_m(x_1,x_2,x_3)t^m$, so that we can write 
\[\prod_{i=1}^3(1-2x_i t)^{-1} = \sum_{m=0}^{\infty} 2^m h_m \cdot t^m\] and \[\prod_{1\leq i<j \leq 3}(1-(x_i+x_j) t)^{-1} = \sum_{m=0}^{\infty} h_m(x_1+x_2,x_1+x_3,x_2+x_3)t^m\] Set $h'_m \coloneqq h_m(x_1+x_2,x_1+x_3,x_2+x_3)$ and $e_m' \coloneqq e_m(x_1+x_2,x_1+x_3,x_2+x_3)$. A direct computation shows that $e_1'=h_1' = 2s_1, e_2' = e_1^2+e_2 = s_2+2s_{1,1},e_3' = e_1e_2-e_3 = s_{2,1}$. 
It holds in general \cite[Theorem 7.6.1]{EC2} that $h_m = e_1h_{m-1}-e_2h_{m-2}+e_3h_{m-3}$ for all $m>0$ so that the same relation holds replacing $h_m,e_m$ with $h_m',e_m'$. So we see that
\begin{align*}
F &= \sum_{m=0}^{\infty} f_m \cdot t^m = \prod_{i=1}^3 \frac{1}{1-2x_i t} \cdot \prod_{1\leq i < j \leq 3} \frac{1}{1-(x_i+x_j)t}\\
&= \left( \sum_{m=0}^{\infty} 2^mh_m t^m \right)\cdot \left(\sum_{m=0} ^{\infty}h'_m t^m \right) \\
&= \left( \sum_{m=0}^{\infty} 2^mh_m t^m \right)\cdot \left( 1+e_1'\sum_{m=1}^{\infty} h'_{m-1} t^m - e'_2\sum_{m=1}^{\infty}h'_{m-2}t^m +e_3'\sum_{m=1}^{\infty}h'_{m-3}t^{m} \right) \\
&= \left( \sum_{m=0}^{\infty} 2^mh_m t^m \right)\cdot \left( 1+e_1't\sum_{m=0} ^{\infty}h'_{m} t^m - e'_2t^2\sum_{m=0}^{\infty}h'_{m}t^m +e_3't^3\sum_{m=0}^{\infty}h'_{m}t^{m} \right) \\
&= \sum_{m=0}^{\infty} 2^m h_m t^m + e_1' F \cdot t - e_2' F \cdot t^2 + e_3' F\cdot t^3 
\end{align*}
and comparing the coefficient in front of $t^m$ on both sides of the above equality, together with the values computed previously, and the fact that $s_m = h_m$ for all $m\geq 0$, concludes the proof. 
\end{proof}

We can now prove Theorem \ref{thm:classV3n}.

\begin{proof}[Proof of Theorem \ref{thm:classV3n}]
By Proposition \ref{prop:V=degeneracy}, $V(3,n)$ is the locus where a morphism $\mathcal{O}^{\oplus n} \to S^2\mathcal{U}^{\vee}$ has rank at most $5$. We can compute the class of this locus via Porteous' formula \cite[Theorem 12.4]{3264} as follows:  Let $\alpha_1,\alpha_2,\alpha_3$ be the Chern roots of $\mathcal{U}^{\vee}$, so that $c_m(\mathcal{U}^{\vee}) = e_m(\alpha) =  e_m(\alpha_1,\alpha_2,\alpha_3)$. By the splitting principle, the Chern roots of $S^2\mathcal{U}^{\vee}$ are $2\alpha_1,2\alpha_2,2\alpha_3,\alpha_1+\alpha_2,\alpha_1+\alpha_3,\alpha_2+\alpha_3$. Denote them as $\beta_1,\dots,\beta_6$ so that $c_m(S^2\mathcal{U}) = e_m(\beta) = e_m(\beta_1,\dots,\beta_6)$. As a consequence of Porteous' formula, the previous discussion, and the dual (or second) Jacobi-Trudi formula, we have  
\[ [V(3,n)] = \det\left(c_{1+j-i}(S^2\mathcal{U}^{\vee})\right)_{1\leq i,j \leq n-5} = \det\left(e_{1+j-i}(\beta)\right)_{1\leq i,j \leq n-5} = h_{n-5}(\beta) \]
Recalling that $\sum_{d=1}^{\infty}h(\beta)t^d = \prod_{i=1}^{6}(1-\beta_it)^{-1}$, we see by definition of the $\beta_i$ that $[V(3,n)] = f_{n-5}(\alpha_1,\alpha_2,\alpha_3)$, where $f_n$ is the symmetric function of Lemma \ref{lem:symmfunctions}. To conclude, we will show that for any $f\in \Lambda_3$ it holds that $f(\alpha_1,\alpha_2,\alpha_3) = [f]$ in $H^{*}(G(3,n))$. It is enough to prove it for $s_1,s_{1,1},s_{1,1,1}$, since they generate $\Lambda_3$ as a ring. This is true because, by \cite[page 178]{3264} the bundle $\mathcal{U}^{\vee}$ has Chern classes $[s_{1^m}] = c_m(\mathcal{U}^{\vee}) = e_m(\alpha_1,\alpha_2,\alpha_3) = s_{1^m}(\alpha_1,\alpha_2,\alpha_3)$. 
\end{proof}

\begin{rem}
With a similar method one can compute the class of the locus $Z_{2}(k,n)$ inside $G(k,n)$. 
\end{rem}

\begin{exa}\label{ex:classesV3n}
Using the formula of Theorem \ref{thm:classV3n}, we can easily write down the classes of the first  interesting cases
\begin{align*}
[V(3,5)]&=[s_0], &[V(3,6)] &= 4[s_1],\\
[V(3,7)] &=11[s_2]+6[s_{1,1}], &[V(3,8)] &= 26[s_3]+23[s_{2,1}]+4[s_{1,1,1}],
\end{align*}
\[[V(3,9)] = 57[s_4]+63[s_{3,1}]+27[s_{2,2}]+18[s_{2,1,1}].\]
Furthermore, a simple Macaulay2 script computes the class of $V(3,n)$ up to $n=100$ in around 2 seconds, and all up to $n=300$ in around 100 seconds on a regular personal computer. 
\end{exa}

Now, we move on to discuss the degree of the projective variety $V(3,n)$. We use the same notation of Theorem \ref{thm:classV3n}.

\begin{cor}\label{cor:degreesV3n}
The degree of $V(3,n)$ viewed as a projective variety under the Plücker embedding is given by the coefficient of $s_{n-3,n-3,n-3}$ in $f_{n-5}\cdot s_1^{2n-4}$.
\end{cor}
\begin{proof}
To obtain the degree of the projective variety $V(3,n)$ under the Plücker embedding, we intersect $V(3,n)$ with $2n-4$ generic Schubert divisors. Thus, the degree is given by the top-degree class $[f_{n-5}\cdot s_1^{2n-4}]\in H^{6(n-3)}(G(3,n))$ under the isomorphism $H^{6(n-3)}(G(3,n))\simeq \mathbb{Z}$ sending $[s_{n-3,n-3,n-3}]\mapsto 1$. \end{proof}

\begin{exa}\label{ex:degreeV3n}
    The first few degrees are given by 
    \begin{align*}
    \deg V(3,5) &=5 ,&\deg V(3,6) &=168 , &\deg V(3,7) &=4032 \\
    \deg V(3,8) &= 84744,&\deg V(3,9) &= 1664091, &\deg V(3,10) &= 31402800
\end{align*}
\end{exa}

Finally, we see that one of the coefficients of the class $[V(3,n)]$ is given explicitly by an Eulerian number:

\begin{cor}\label{cor:euler}
    If $n\geq 5$, the coefficient of $[s_{n-5}]$ in the class $[V(3,n)]$ is the Eulerian number $A(n-3,1)=2^{n-3}-(n-2)$ which  counts the number of permutations in $\mathfrak{S}_n$ with $1$ descent.
\end{cor}
\begin{proof} This can be verified directly for $5\leq n \leq 9$ thanks to  Example \ref{ex:classesV3n}. In general, we prove the statement by induction on $n$, using Theorem \ref{thm:classV3n}. Let $A_n$ be the coefficient of $[s_{n-5}]$ in the class $[V(3,n)] = [f_{n-5}]$. This is also the coefficient of $s_{n-5}$ in the Schur function expansion of $f_{n-5}$. Let $n\geq 9$ and  and assume that the statement holds for all values smaller than $n$. The recursive formula of Theorem \ref{thm:classV3n}, together with Pieri's formula shows that $A_n=2A_{n-1}-A_{n-2}+2^{n-5}$. Induction shows that $A_n = 2(2^{n-4}-(n-3))-(2^{n-4}-(n-4))+2^{n-5} = 2^{n-3}-(n-2)$. 
\end{proof}
\begin{rem}\label{rem:physics}
This last corollary can be interpreted and generalized by the formalism of spinor-helicity in physics, following \cite{lam2024moduli} and \cite{EMPS24}. Fix a bilinear pairing on $V = \mathbb{C}^{1\times n}$, so that we have an identification $V\cong V^{\vee}$ and for each subspace $U\subseteq V$ we can consider the annihilator $U^{\perp}\subseteq V$. The spinor-helicity variety is
\[ SH(2,n,0) = \left\{ (U_1,U_2) \in  G(2,n)\times G(2,n) \,|\, U_1\subseteq U_2^{\perp}  \right\}. 
\]
For each $0\leq d \leq n-4$ the lifted scattering correspondence is 
\[ \widetilde{C}_{d+1} = \{ ((U_1,U_2),W) \in S(2,n,0)\times V(d+1,n) \,|\, U_1 \subseteq W \subseteq U_2^{\perp} \}.  \]
The map $\widetilde{C}_{d+1} \to SH(2,n,0)$ is generically finite, and the fibers correspond to solutions to the Cachazo-He-Yuan scattering equations \cite{CHY13}. In particular, Cachazo-He-Yuan \cite{CHY13} proved that the degree of the map is the Eulerian number $A(n-3,d-1)$, see \cite[Theorem 6.3]{EMPS24} for a proof. This can be interpreted as follows: take a general full flag $V_{\bullet} = 0\subseteq V_1 \subseteq \dots \subseteq V_n=V$ and consider the Schubert variety $\Sigma_{((n-d-1)^2,2^{d-1})}(V_{\bullet}) = \Sigma_{(n-d-1,n-d-1,2,\dots,2)}(V_{\bullet})$ of this flag. Setting $U_1=V_2$ and $U_2 = V_{n-2}^{\perp}$, we see that $(U_1,U_2)$ is a general point in $S(2,n,0)$, and moreover
\[ V(d+1,n)\cap \Sigma_{(n-d-1)^2,2^{d-1}}(V_{\bullet}) = \{W\in V(d+1,n) \,|\, U_1 \subseteq W \subseteq U_2^{\perp}\}.  \]
This set is precisely a general fiber of the map $\widetilde{C}_{d+1}\to SH(2,n,0)$ and the results mentioned before show that it consists of $A(n-3,d-1)$ points. As a consequence, we get:
\end{rem}

\begin{prop}\label{prop:eulergeneral}
The coefficient of $[s_{(n-d-3)^{d-1}}]$ in the class of $[V(d+1,n)]$ in $G(d+1,n)$ is $A(n-3,d-1)$.
\end{prop}
\begin{proof}
By \cite[Proposition 4.6]{3264}, this coefficient is given by the number of intersection points of $V(d+1,n)$ with a general Schubert variety $\Sigma_{(n-d-1)^2,2^{d-1}}$. Then the result follows from Remark \ref{rem:physics}.
\end{proof}

\begin{rem}
We conclude this section observing that any result obtained for $V(d+1,n)$ can be translated into the analogous result for $V(n-d-1,n)$. Indeed, fix the standard bilinear form on $V=\mathbb{C}^{1\times n }$ and the corresponding duality 
\[ G(d+1,n) \to G(n-d-1,n), \qquad [U] \mapsto [U^{\perp}]. \]
Then this induces an isomorphism $V(d+1,n) \cong V(n-d-1,n)$. This is the content of Goppa's duality, stating that if $n$ points lie on a rational normal curve in $\mathbb{P}^{d}$, then their Gale transform consists of $n$ points lying on a rational normal curve in $\mathbb{P}^{n-d-2}$. For a reference, see \cite[Corollary 3.2]{EP00} and the remarks afterward.
\end{rem}

\section{A stratification of $V(k,n)$}\label{sec:strat}
The Grassmannian has a matroidal stratification defined as follows: for a point $W\in G(d+1,n)$, its matroid, considered as a set of bases, is defined as 
\begin{gather*}
    \mathcal M(W) \coloneqq \left\{ I \in \binom{[n]}{d+1}\, |\, p_I(W) \neq 0 \right\}
\end{gather*}
where $p_{I}(W)$ denotes the Pl\"ucker coordinate indexed by $I$. If instead all the Pl\"ucker coordinates of $W$ are real and non-negative, then $W$ belongs to the non-negative Grassmannian $G_{\geq 0}(d+1,n)$ and $\mathcal{M}(W)$ is called a positroid. Given a matroid (resp. positroid) of the form $\mathcal{M} = \mathcal{M}(W)$, the corresponding matroid (resp. positroid) cell is $\Pi_{\mathcal{M}}=\{X\in G(d+1,n) \,|\, \mathcal{M}(X)=\mathcal{M}\}$. These cells give a stratification of $G(d+1,n)$ and $G_{\geq 0}(d+1,n)$ and their closure is
$\overline{\Pi}_{\mathcal{M}} = \{X\in G(d+1,n) \,|\, p_I(X) = 0\, \text{ if } I\notin \mathcal{M}\}$.
We refer the reader to \cite{lam2015totally} for a rigorous exposition on the topic. 

\begin{rem}\label{rem:pluecker}
    Consider now the Veronese map $\theta_{d}\colon G(2,n) \dashrightarrow G(d+1,n)$ and let $W\in G(2,n)$ be such that the map is defined. Let also $1\leq i_1<\dots <i_{d+1} \leq n$ be a subset with $d+1$ elements: the formula for the Vandermonde determinant shows 
\begin{equation}\label{eq:vandermonde} 
    p_{i_1\dots i_{d+1}}(\theta_{d}(W)) = \prod_{1\leq a < b \leq d+1} p_{i_ai_b}(W)  
    \end{equation}
    In terms of the corresponding matroids this means that:
    \[ \mathcal{M}(\theta_{d}(W)) = \left\{ I \in \binom{[n]}{d+1} \,|\, \binom{I}{2} \subseteq \mathcal{M}(W)  \right\}\]
    This leads us to the following definition:
\end{rem}

\begin{dfn}[Veronese image of a rank $2$ matroid]
Let $\mathcal{M}$ be a matroid of rank $2$ on $[n]$ and let $k\geq 2$. Its $d$-th Veronese image is
\[ \theta_{d}(\mathcal{M}) = \left\{ I \in \binom{[n]}{k} \,|\, \binom{I}{2} \in \mathcal{M}  \right\} \]
\end{dfn}

\begin{exa}
When $d = 2$ we have 
$$\theta_{2}(\mathcal M) = \{ \{i,j,k\} \subset [n] \; |\; \{i,j\}, \{j,k\}, \{i,k\} \in \mathcal M \}.$$  
\end{exa}

\begin{prop}
    Given a matroid $\mathcal M$ of rank $2$ on $[n]$ and $k
    \geq 2$, its $d$-th Veronese image $\theta_{d}(\mathcal{M})$ is either empty or a matroid. 
\end{prop}
\begin{proof}
    Assume that $\theta_d(\mathcal{M})$ is not empty. We want to show that $\theta_d(\mathcal{M})$ is in fact a matroid by showing that the exchange axiom is satisfied, which is true if $\theta_d(\mathcal M)$ contains a single element. So, we assume there are at least two elements in $\theta_d(\mathcal M)$. 
    
    Let $I, J$ be $d+1$ subsets of $[n]$ contained in $\theta_d(\mathcal{M})$. Let $a \in I \setminus J$. We want to show that there exists some $b \in J\setminus I$ such that $J \setminus \{b\} \bigcup \{a\} \in \theta_d (\mathcal{M})$, meaning that $\{a,j\} \in \mathcal{M}$ for all $j\in J\setminus\{b\}$.
    Let $a'\in I$ be distinct from $a$ and consider $\{a,a'\}\in \mathcal{M}$. This needs to satisfy the exchange axiom with any $2$-subset $\{j_1,j_2\}$ of $J$, hence $\{a,j_1\}\in \mathcal{M}$ or $\{a,j_2\}\in \mathcal{M}$. Thus there can be at most one element $b\in J$ such that $\{a,b\}\notin \mathcal{M}$. If such an element $b$ exists, then $J\setminus \{b\} \cup \{a\} \in \mathcal{M}$ by the previous discussion. If there is no such element, take any $b\in J\setminus I$.
\end{proof}        

    In addition to being a matroid, $\theta_{d}(\mathcal M(X))$ is a positroid if $X \in G_{\geq 0}(2,n)$. This follows from Remark \ref{rem:pluecker}, since the product of positive $2 \times 2$ minors is also positive. 
Thus, for each matroid (resp. positroid) $\mathcal{M}$ coming from $G(2,n)$ we obtain another matroid (resp. positroid) $\theta_{d}(\mathcal{M})$ on $G(d+1,n)$ and the corresponding strata $\Pi_{\theta_{d}(\mathcal{M})}$ on $G(d+1,n)$. In particular, Remark \ref{rem:pluecker} shows that
\[ \theta_{d}(\Pi_{\mathcal{M}}) \subseteq \Pi_{\theta_{d}(\mathcal{M})} \cap V(d+1,n) \]
We can use this to investigate whether $V(d+1,n)$ is a positive geometry. More precisely Thomas Lam conjectured in \cite[Conjecture 4.10]{lam2024moduli} that the closure 
\[ V_{\geq 0}(d+1,n) = \overline{\theta_{d}(G_{>0}(2,n))}\] 
should be a positive geometry. In this case, it is natural to expect, according to Heuristic 4.1 in \cite{arkani2017positive} , that the map
\[ \theta_{d}\colon (G(2,n),G_{\geq 0}(2,n)) \dashrightarrow (V(d+1,n),V_{\geq 0}(d+1,n)) \]
should be a morphism of positive geometries. In particular the boundary of $G_{\geq 0}(2,n)$ should be mapped to the boundary of $V_{\geq 0}(d+1,n)$. Experiments with computer algebra suggest however, that the boundary of $G_{\geq 0}(2,n)$ is not enough to obtain the whole boundary of $V_{\geq 0}(d+1,n)$. More precisely, we consider $d=2$ and $n=6$. One 6-dimensional component of the boundary of $G_{\geq 0}(2,6)$ is given by the closure of $\Pi_{\mathcal{M}_1}$ where $\mathcal{M}_1= \{\{i,j\}\subset [6] \,|\, i,j\ne 1\}$, and its $5$-dimensional boundary is given by closures of strata $\Pi_{\mathcal{M}_2}$, where $\mathcal{M}_2$ has the form $\mathcal{M}_2= \mathcal{M}_1 \setminus \{\{i,j\}\}$, for any specific choice of $\{i,j\} \in \mathcal{M}_1$. A computation with OSCAR \cite{OSCAR}, shows that 
\begin{align*}
\dim \overline{\theta_{2}(\Pi_{\mathcal{M}_1})} &= \dim \overline{\Pi_{\theta_2(\mathcal{M}_1)}} \cap V(3,6) =  6, \\  
\dim \overline{\theta_{2}(\Pi_{\mathcal{M}_2})} &= \dim \overline{\Pi_{\theta_2(\mathcal{M}_2)}} \cap V(3,6) =  4.\end{align*}
Hence, even if $\Pi_{\mathcal{M}_2}$ gives a codimension one boundary of $\Pi_{\mathcal{M}_1}$, this is not true anymore under the Veronese mapping. Thus, if $V(d+1,n)$ is a positive geometry, there must be other components in the boundary of $\overline{\theta_{d}(\Pi_{\mathcal{M}_1})}$, not coming from $G_{\geq 0}(2,n)$. This is likely due to the fact that the Veronese map is not defined everywhere: to get an actual map, we would need to blow up the indeterminacy locus, and then the extra components that we are missing should probably arise from the exceptional divisor.



\section*{Acknowledgements}

We first came in contact with the ABCT varieties in Thomas Lam's course at the School \textit{Combinatorial Algebraic Geometry from Physics}, held at the MPI-MiS Leipzig in 2024. We wish to thank Thomas Lam and the organizers of the School. We also wish to thank Thomas Lam and Bernd Sturmfels for encouragement and advice. We thank Fulvio Gesmundo, Leonie Kayser,  M\'at\'e Telek and Emanuele Ventura for their interest and helpful comments.



\bibliography{bibliography}

\noindent

\nocite{*}
\bibhere

\end{document}